\documentclass[a4paper,12pt,leqno]{article}
\usepackage{latexsym}
\usepackage[all]{xy}

\usepackage{amssymb} 
\usepackage{amsmath} 
\usepackage{theorem}

\def\Z{{\mathbb{Z}}}
\def\C{{\mathbb{C}}}
\def\R{{\mathbb{R}}}
\def\K{{\mathbb{K}}}
\def\A{{\mathcal{A}}}

\DeclareMathOperator{\Der}{Der}

\numberwithin{equation}{section}

\newcommand{\owari}{\hfill$\square$}

\theoremstyle{break}
\newtheorem{theorem}{Theorem}[section]
\newtheorem{prop}[theorem]{Proposition}
\newtheorem{cor}[theorem]{Corollary}
\newtheorem{lemma}[theorem]{Lemma}
\newtheorem{define}[theorem]{Definition}
\newtheorem{rem}[theorem]{Remark}

\newtheorem{example}[theorem]{Example}

\newcommand{\xgraphAvertex}[1][****]{
\xgraphAVertex #1
}
\newcommand{\xgraphAVertex}[4]{
\if#3o\put(0,0){\circle{4}}\fi%
\if#3*\put(0,0){\circle*{4}}\fi%
\if#3.\put(0,0){\circle*{2.4}}\fi%
\if#4o\put(30,0){\circle{4}}\fi%
\if#4*\put(30,0){\circle*{4}}\fi%
\if#4.\put(30,0){\circle*{2.4}}\fi%
\if#2o\put(0,30){\circle{4}}\fi%
\if#2*\put(0,30){\circle*{4}}\fi%
\if#2.\put(0,30){\circle*{2.4}}\fi%
\if#1o\put(30,30){\circle{4}}\fi%
\if#1*\put(30,30){\circle*{4}}\fi%
\if#1.\put(30,30){\circle*{2.4}}\fi%
}
\newcommand{\xgraphA}[6]{
\if#1+\put(30,30){\line(-1,0){30}}\fi
\if#1.\qbezier[7](30,30)(15,30)(0,30)\fi
\if#1-\put(30,31){\line(-1,0){30}}\put(30,29){\line(-1,0){30}}\fi
\if#2+\put(0,0){\line(1,1){30}}\fi 
\if#2.\qbezier[10](0,0)(15,15)(30,30)\fi 
\if#2-\put(-0.7,0.7){\line(1,1){30}}\put(0.7,-0.7){\line(1,1){30}}\fi 
\if#3+\put(30,30){\line(0,-1){30}}\fi
\if#3.\qbezier[7](30,0)(30,15)(30,30)\fi 
\if#3-\put(29,30){\line(0,-1){30}}\put(31,30){\line(0,-1){30}}\fi
\if#4+\put(0,0){\line(0,1){30}}\fi 
\if#4.\qbezier[7](0,0)(0,15)(0,30)\fi 
\if#4-\put(-1,0){\line(0,1){30}}\put(1,0){\line(0,1){30}}\fi 
\if#5+\put(0,30){\line(1,-1){30}}\fi 
\if#5.\qbezier[10](30,0)(15,15)(0,30)\fi 
\if#5-\put(-0.7,29.3){\line(1,-1){30}}\put(0.7,30.7){\line(1,-1){30}}\fi 
\if#6+\put(0,0){\line(1,0){30}}\fi 
\if#6.\qbezier[7](0,0)(15,0)(30,0)\fi 
\if#6-\put(0,1){\line(1,0){30}}\put(0,-1){\line(1,0){30}}\fi 
\xgraphAvertex}

\title{A generalized logarithmic module and duality of Coxeter multiarrangements}
\author{Takuro Abe\thanks{
Department of Mathematics, Hokkaido University, 
Kita-10, Nishi-8, Kita-Ku, 
Sapporo, Hokkaido 060-0810, Japan.
email:abetaku@math.sci.hokudai.ac.jp.}}
\date{\today}

\begin{document}

\maketitle

\begin{abstract}
We introduce a new definition of a generalized logarithmic module of 
multiarrangements by uniting those of 
the logarithmic derivation and the differential modules. 
This module is realized as a logarithmic derivation module of an arrangement of 
hyperplanes with a multiplicity consisting of both positive and negative integers. We  
consider 
several properties of this module including Saito's criterion and reflexivity. 
As applications, we prove a shift isomorphism and 
duality of some Coxeter multiarrangements by using the primitive 
derivation. 
\end{abstract}
\setcounter{section}{-1}

\section{Introduction}
Let $V$ be an $\ell$-dimensional vector space over the real number field $\R$,  
$\{x_1,\ldots,x_\ell\}$ a basis for the dual vector space $V^*$, 
and 
$S:=\mbox{Sym}(V^*) \otimes_\R \C
 \simeq \C[x_1,\ldots,x_\ell]$. 
Let $\Der_{\C}(S)$ denote 
the $S$-module of $\C$-linear derivations of $S$ and 
$\Omega^1_V$ the $S$-module of differential $1$-forms, 
i.e., 
$
\Der_{\C}(S)=\bigoplus_{i=1}^{\ell} S \cdot \partial_{x_i}$ 
and 
$\Omega^1_V:=\bigoplus_{i=1}^\ell S\cdot dx_i$. 
A non-zero element $\theta=\sum_{i=1}^\ell f_i \partial_{x_i} 
\in \Der_{\C}(S)$ (resp. $\omega=\sum_{i=1}^\ell g_i d{x_i} 
\in \Omega^1_V)$ is 
\textit{homogeneous of degree $p$} if $f_i\ (\mbox{resp}.\ g_i)$ is zero or homogeneous of 
degree $p$ for each $i$. 

A \textit{hyperplane arrangement} $\A$ (or simply 
an \textit{arrangement}) is a finite collection of affine hyperplanes 
in $V$. If each hyperplane in $\A$ contains the origin, we say that $\A$ is 
\textit{central}. In this article we assume that all arrangements are central 
unless otherwise specified. 
A \textit{multiplicity} $m$ on an arrangement $\A$ is a map 
$m:\A \rightarrow \Z_{\ge0}$ and a pair $(\A,m)$ is called a 
\textit{multiarrangement}. 
Let $|m|$ denote the sum of the multiplicities 
$\sum_{H \in \A}m(H)$. When $m \equiv 1,\ (\A,m)$ is the same as the hyperplane arrangement 
$\A$ and sometimes called a \textit{simple arrangement}. For each hyperplane $H \in \A$ 
fix a linear form $\alpha_H \in V^*$ such that $\ker(\alpha_H)=H$. Put 
$
Q(\A,m):=\prod_{H \in \A} \alpha_H^{m(H)}.
$
The main objects in this article are the \textit{logarithmic derivation module} 
$D(\A,m)$ of $(\A,m)$ defined by 
$$
D(\A,m):=\{\theta \in \Der_{\C}(S)|\theta(\alpha_H) \in S\cdot \alpha_H^{m(H)}\ 
\mbox{for all }H \in \A \},
$$
and the \textit{logarithmic differential module} $\Omega^1(\A,m)$ of $(\A,m)$ defined by 
$$
\Omega^1(\A,m):=\{\omega \in \displaystyle \frac{1}{Q(\A,m)} \Omega^1_V|d \alpha_H \wedge \omega\ 
\mbox{is regular along }H\ \mbox{for all }H \in \A\}.
$$
It is well-known that $D(\A,m)$ and $\Omega^1(\A,m)$ are $S$-dual modules, and hence 
reflexive in general. 
A multiarrangement $(\A,m)$ is \textit{free} if 
$D(\A,m)$ is a free $S$-module of rank $\ell$. If $(\A,m)$ is free, then 
there exists a homogeneous free basis $\{\theta_1,\ldots, \theta_\ell\}$ for 
$D(\A,m)$. Then we define the \textit{exponents} of a free multiarrangement 
$(\A,m)$ by $\exp(\A,m):=(\deg(\theta_1),\ldots,\deg(\theta_\ell))$. 
The exponents are independent of a 
choice of a basis. When $m \equiv 1$, 
the logarithmic derivation, differential modules and exponents are denoted by 
$D(\A)$, $\Omega^1(\A)$ and $\exp(\A)$ respectively. 

Multiarrangements were introduced by Ziegler in \cite{Z}, and 
have been shown to illuminate several algebraic and topological problems about 
simple arrangements. 
Two important results concerning multiarrangements are the freeness of 
Coxeter arrangements with (quasi-)constant multiplicities (\cite{ST2}, \cite{T4}, \cite{Y0} and 
\cite{AY2}), and the relation with the Hodge filtration (\cite{T5}). 
In these results, only multiplicities of positive, or non-negative integers were considered. 
However, the results 
in \cite{AY2} and \cite{ANN} suggest (at least when we consider 
Coxeter arrangements) mixing 
the definitions of the logarithmic derivation and the differential modules 
by considering a multiplicity with positive and negative integers. 
Let us explain the background more explicitly. 

Let $\A$ be a Coxeter arrangement 
with the Coxeter number $h$ and $m:\A \rightarrow \{0,1\}$ a 
\textit{quasi-constant} multiplicity. Then, 
for a multiarrangement $(\A,2k \pm m)$ with $k \in \Z_{> 0}$, 
Theorem 10 in 
\cite{AY2} gives isomorphisms 
\begin{eqnarray*}
D(\A,m)(-kh) \rightarrow D(\A,2k+m)\\
\Omega^1(\A,m)(-kh) \rightarrow D(\A,2k-m)
\end{eqnarray*}
and a duality between $D(\A,2k+m)(kh)$ and $D(\A,2k-m)(kh)$
by using Kyoji Saito's primitive derivation $D$ 
and the affine connection $\nabla$, 
where, for an $S$-module $M$, $M(a)$ is an $S$-graded module such that 
its $b$-degree part $M(a)_b=M_{a+b}$ for $a,b \in \Z$. 
For details of $D$ and $\nabla$, see \cite{Sa1} and Section two. 
Also, the results in \cite{ANN} shows that, 
for a braid arrangement $\A$ and multiplicity $m:\A \rightarrow \{+1,0,-1\}$, there exists 
a similar duality 
between free multiarrangements $(\A,2k + m)$ and $(\A,2k - m)$. 
Hence it is natural to expect a similar result to that in \cite{AY2} 
for a multiplicity 
$m:\A \rightarrow \{+1,0,-1\}$ on an arbitrary Coxeter arrangement $\A$.
For that purpose, we need a generalized logarithmic module for 
such a 
multiplicity.
Our first goal of this article is to define 
this module for a multiplicity with both positive and negative integer values. 
Let us begin with a generalized definition of multiplicities. 

\begin{define}
We call a map $m:\A \rightarrow \Z$ 
a \textit{multiplicity} on $\A$, and 
$(\A,m)$ a \textit{multiarrangement}. Define 
$\A_{+}:=\{H \in \A|m(H) >0\}$, 
$\A_{-}:=\{H \in \A|m(H) <0\}$, and two polynomials $Q_+$ and $Q_-$ by 
\begin{eqnarray*}
Q_+:  &=&\prod_{H_+ \in \A_+}\alpha_{H_+}^{m(H_+)},\\
Q_{-}:&=&\prod_{H_- \in \A_{-}}\alpha_{H_-}^{-m(H_-)}.
\end{eqnarray*}

For $\mu \in \{+,-\}$, define $m_\mu$ by 

\[
m_\mu(H):=
\left\{
\begin{array}{rl}
m(H) & \mbox{if}\ H \in \A_\mu ,\\
0\ \     & \mbox{if}\ H \not \in \A_\mu.
\end{array}
\right.
\]

We say that $\displaystyle \frac{Q_+}{Q_-}$ is a ``defining function'' of 
a multiarrangement $(\A,m)$.
\end{define}

In order to define a logarithmic module, 
let us introduce some morphisms between localizations of $\Der_\C(S)$ and 
$\Omega^1_V$ (see \cite{Sa1}, \cite{Sa3} and \cite{Sa4} for details). 
If we fix an inner product 
$I^*$ of $V^*$ then the dual isomorphism $I^*:\Omega^1_V \rightarrow \Der_\C(S)$ is 
canonically defined by $I^*(\omega)(f):=
I^*(\omega, df )$ for $f \in S$. Note that $\Omega^1(\A_-,-m_-) \subset (\Omega^1_V)_{(0)}$ and 
$D(\A_+,m_{+}) \subset \Der_\C(S) \subset \Der_\C(S)_{(0)}$. 
If we agree that $I^*$ also denotes the localized identification 
$(\Omega^1_V)_{(0)} \rightarrow \Der_\C(S)_{(0)}$ 
at the 
ideal $(0)$, then we obtain the following diagram:

$$\xymatrix@R1pc{
\Omega^1(\A_-,-m_-)\  \ar @{^{(}-{>}}[d] & \displaystyle \frac{1}{Q_-} D(\A_+,m_+) \ar @{^{(}-{>}}[d] \\
(\Omega^1_V)_{(0)} \ar[r]^{I^*} & \Der_\C(S)_{(0)}  }$$

Now we can state a definition of the main character of this article.

\begin{define}
The \textit{(generalized) logarithmic module} $D\Omega(\A,m)$ of a multiarrangement $(\A,m)$ is defined, 
as the submodule of $\Der_\C(S)_{(0)}$, by 
$$
D\Omega(\A,m):=\displaystyle \frac{D(\A_+,m_{+})}{Q_{-}} \bigcap I^*(\Omega^1(\A_{-},-m_{-})).
$$
\label{DO}
\end{define}

\begin{rem}
Since $I^*$ is a canonical identification, we can also define an $S$-module 
$\Omega D(\A,m)$, as the submodule of $(\Omega^1_V)_{(0)}$, by 
$$
\Omega D(\A,m):=(I^*)^{-1}
(\displaystyle \frac{D(\A_-,-m_{-})}{Q_{+}}) \bigcap \Omega^1(\A_{+},m_{+}).
$$
Since $D \Omega (\A,m) \simeq \Omega D(\A,-m)$,
in this article, we 
mainly study and use $D \Omega (\A,m) \subset \Der_\C(S)_{(0)}$.
\end{rem}

\begin{rem}
We defined the generalized logarithmic module over the 
real number field and its complexification, but we can give that 
definition over any field $\K$ by introducing 
a non-degenerate symmetric bilinear form $I^*:\K^\ell \times \K^\ell \rightarrow \K$. 
In particular, 
if we fix an orthonormal coordinate system $\{x_1,\ldots,x_\ell\}$ for $V^*$, 
Definition \ref{DO} is equivalent to 
the following: 

\begin{define}
The logarithmic module $D\Omega(\A,m) \subset 
\displaystyle \frac{1}{Q_-} \Der_{\K}(S)$ is defined by 
$
D\Omega(\A,m):=\{\theta \in \displaystyle \frac{1}{Q_-}\Der_{\K}(S)\ |
\ \theta(\alpha_{H_+} ) \in \displaystyle \frac{\alpha_{H_+}^{m(H_+)}}{Q_-} \cdot S\ 
\mbox{for all}\ H_+ \in \A_+\ \mbox{and}\ 
\theta \wedge \partial \alpha_{H_-}\ \mbox{is regular}\ \\
\mbox{along}\ 
H_- \ \mbox{for all}\ H_- \in \A_-\}$, where 
$\partial \alpha_H:=\sum_{i=1}^\ell \partial_{x_i}(\alpha_H)\partial_{x_i}$. 
\label{DO2}
\end{define}

Note that, here, we are canonically identifying $\partial_{x_i}$ and $dx_i$. 
Depending on a choice of (non-orthonormal) coordinates, the defining 
statement of Definition \ref{DO2} 
changes, though Definition \ref{DO} does not.
For it includes 
inner products in the definition. In this article, we work over the 
real number field and its complexification, 
but almost all results (more precisely, results except for those in Section two) 
hold true even over any fields by applying the same proofs as for those over $\R$ 
or $\C$. 

Also, when we prove results related to $D\Omega(\A,m)$ we 
often choose an orthonormal coordinate system 
$\{x_1,\ldots,x_\ell\}$ and use Definition \ref{DO2} 
for convenience. 
\end{rem}

\begin{rem}
The $S$-module structure of $D\Omega(\A,m)$ depends on the choice of the inner product $I^*$, which 
we will see in Example \ref{A1A1}. Hence we sometimes write $D\Omega(\A,m,I^*)$ instead of 
$D\Omega(\A,m)$ when we want to make it clear with which inner products we are studying. When there 
would be no confusions, we omit $I^*$ and just write $D\Omega(\A,m)$. 
\label{inner}
\end{rem}

The aim of this article is to investigate properties of the logarithmic module $D\Omega(\A,m)$ 
(Theorems \ref{Saito} and \ref{duality} for instance), 
and give an application to Coxeter multiarrangements as the main theorem and 
a generalization of results in both \cite{AY2} and \cite{ANN}. 
More explicitly, we will give 
the following shift isomorphism by using the primitive derivation:

\begin{theorem}
Let $\A$ be a Coxeter arrangement with the Coxeter number 
$h$, $m:\A \rightarrow \{+1,0,-1\}$ and $k \in \Z_{>0}$. 
Then there are $S$-module isomorphisms as follows:
\begin{eqnarray*}
D\Omega(\A,m)(-kh) &\rightarrow& D(\A,2k+m),\\
\Omega D(\A,m)(-kh) &\rightarrow& D(\A,2k-m).
\end{eqnarray*}
\label{mainmain}
\end{theorem}

We prove Theorem \ref{mainmain} by constructing explicit isomorphisms in Theorem \ref{Coxeter}. 
In particular, 
the generalized duality of Coxeter multiarrangements 
between $D(\A,2k+m)(kh)$ and $D(\A,2k-m)(kh)$ 
will be given in Corollary \ref{coxeterduality}. 

The organization of this article is as follows. In Section one we show some general properties of 
$D\Omega(\A,m)$ including Saito's criterion and reflexivity. In Section two we apply the theory of 
logarithmic modules to Coxeter multiarrangements and prove Theorem \ref{mainmain} and 
the duality. 
\medskip

\noindent
\textbf{Acknowledgements}. The author is grateful to Professor Hiroaki Terao and 
Professor Masahiko Yoshinaga for 
useful discussions and suggestions to this article. Also the author thanks 
Professor Sergey Yuzvinsky and Dr. Max Wakefield for their kind advice. 
The author is supported by the 
JSPS Research Fellowship for Young Scientists.

\section{Generalized logarithmic module of a multiarrangement}
In this section we study several properties of a logarithmic module $D\Omega(\A,m)$ which is 
introduced in the previous section. Let $V$ be an $\ell$-dimensional Euclidean space. 
Before fixing an inner product, 
let us consider the following example.


\begin{example}
Let $(\A,m)$ be a multiarrangement in $\R^2$ defined by 
$\displaystyle \frac{x}{y}=0$, where $\{x,y\}$ 
is a basis  for $(\R^2)^*$. Let $I_1^*$ be an 
inner product on $V^*$ defined by 
$$
I_1^*(dx,dx)=I_1^*(dy,dy)=1,\ I_1^*(dx,dy)=I_1^*(dy,dx)=0.
$$
Then it is easy to show that $D\Omega(\A,m,I_1^*)$ is a free $S$-module with basis 
$$
D\Omega(\A,m,I_1^*)=\langle x \partial_x, \displaystyle \frac{1}{y} \partial_y \rangle_S.
$$
Next, let $I_2^*$ be an inner product on $(\R^2)^*$ 
defined by 
$$
I_2^*(dx,dx)=I_2^*(dx,dy)=I_2^*(dy,dx)=1,\ I_2^*(dy,dy)=2.
$$
Then it is also easy to show that $D\Omega(\A,m,I_2^*)$ is a free $S$-module with basis 
$$
D\Omega(\A,m,I_2^*)=\langle \partial_y, \displaystyle \frac{x \partial_x + 2 x \partial_y}{y}
  \rangle_S = \langle I_2^*(-dx+dy), I^*_2(\displaystyle \frac{x}{y} dy) \rangle_S.
$$
Hence we can see that 
the $S$-module structure of $D\Omega(\A,m,I^*)$ depends on a choice of 
inner product $I^*$ as mentioned in 
Remark \ref{inner}. 
\label{A1A1}
\end{example}
 
From now on, in this section, let us fix an inner product $I^*$ on $V^*$ and 
write $D\Omega (\A,m)$ instead of $D\Omega (\A,m,I^*)$ unless otherwise specified. 
First, note that, by definition $D\Omega(\A,m)$ is an $S$-graded module. 
If $Q_+=1$ then Definition \ref{DO} coincides with that of $\Omega^1(\A,m)$ 
and 
if $Q_-=1$ then Definition \ref{DO} coincides with that of $D(\A,m)$ by the 
fixed identification between $\Der_\C(S)$ and $\Omega^1_V$ through $I^*$. 
Hence $D\Omega(\A,m)$ is a generalization of both modules 
$D(\A,m)$ and $\Omega^1(\A,m)$.

\begin{example}
Let $(\A,m)$ be a multiarrangement in $\R^2$ defined by 
$\displaystyle \frac{y}{x(x-y)}=0$, where $\{x,y\}$ 
is an orthonormal coordinate system for $\R^2$. Then it is easy to prove that 
$D\Omega(\A,m)$ is free with a basis 
$$
D\Omega(\A,m)=\langle \displaystyle \frac{\partial_x}{x},\displaystyle \frac{x \partial_x - y \partial_y}{x-y} \rangle_S.
$$
\label{B2}
\end{example}

Though two examples above are both free (see Corollary \ref{2multi}), 
it is known that $D\Omega(\A,m)$ is not free in general when $\ell \ge 3$ 
(see Example \ref{A3} for example). Hence we say that 
$(\A,m)$ is \textit{free} if $D\Omega(\A,m)$ is a free $S$-module of rank $\ell$. 
If $(\A,m)$ is free then there exists a homogeneous free basis $\{\theta_1,\ldots,\theta_\ell\}$ for 
$D\Omega(\A,m)$.
Then the \textit{exponents} of a free multiarrangement 
$(\A,m)$ is defined by $\exp(\A,m):=(\deg \theta_1,\ldots,\deg \theta_\ell)$. Exponents are independent of 
a choice of a basis. 
For instance, the multiarrangement in Example 
\ref{B2} is free with exponents $(-1,0)$. 

Next let us consider the structure of the module $D\Omega(\A,m)$, especially its freeness.  


\begin{lemma}
For $\theta_1,\ldots,\theta_\ell \in D\Omega(\A,m)$, let us define the 
$(\ell \times \ell)$-matrix $M(\theta_1,\ldots,\theta_\ell)$ as 
its $(i,j)$th entry $\theta_j(x_i)$. Then 
$$
\det M(\theta_1,\ldots,\theta_\ell) \in S \cdot \displaystyle \frac{Q_+}{Q_-}.
$$
\label{wedge}
\end{lemma}

\noindent
\textbf{Proof}. We prove by choosing an orthonormal basis $\{x_1,\ldots,x_\ell\}$ for 
$V^*$. 
By definition, 
poles of $\det M(\theta_1,\ldots,\theta_\ell)$ can exist only along 
$Q_-=0$. 
Take $H_0 \in \A$ and put  
$|m(H_0)|=m_0$. First, assume that $Q_+=Q_+'\alpha_{H_0}^{m_0}$ with 
$ \alpha_{H_0} \nmid Q_+$. 
By an appropriate change of orthonormal coordinates, 
we may assume that $\alpha_{H_0}=x_1$.
Then $\theta_i(x_1) \in \displaystyle \frac{x_1^{m_0}}{Q_-} \cdot S$ 
for $i=1,2,\ldots,\ell$. 
Therefore, 
$\det M(\theta_1,\ldots,\theta_\ell) \in S \cdot \displaystyle \frac{x_1^{m_0}}{Q_-^\ell}$. 
Since $H_0$ runs over all $H_0 \in \A_+$, it holds that 
$\det M(\theta_1,\ldots,\theta_\ell) \in S \cdot \displaystyle \frac{Q_+}{Q_-^\ell}$.
Next, assume that 
$Q_-=\alpha_{H_0}^{m_0} Q_-'$ with $\alpha_{H_0} \nmid Q_-'$. Again we may assume that 
$\alpha_{H_0}=x_1$. 
Let us show that the order of the pole of $\det M(\theta_1,\ldots,\theta_\ell)$ along $H_0$ 
is at most $m_0$. By definition $\theta_j \wedge d x_1$ has no poles along $H_0$ for all $j$. 
In particular, $\theta_j(x_i)$ is regular along $H_0$ for $j=1,2,\ldots,\ell$ and 
$i=2,\ldots,\ell$. 
Hence 
the order of poles of $\det M(\theta_1,\ldots,\theta_\ell)$ along $H_0$ is at most $m_0$, 
which shows that 
$\det M(\theta_1,\ldots,\theta_\ell) \in S \cdot \displaystyle \frac{Q_+}{Q_-}$. 
\owari
\medskip

Now we obtain Saito's criterion for  
$D\Omega(\A,m)$.

\begin{theorem}[Saito's criterion]
For $\theta_1,\ldots,\theta_\ell \in D\Omega(\A,m)$, the following two conditions are equivalent:
\begin{itemize}
\item[(1)] $\{\theta_1,\ldots,\theta_\ell\}$ forms a basis for $D\Omega(\A,m)$.
\item[(2)] $\det M(\theta_1,\ldots,\theta_\ell)=c \cdot 
\displaystyle \frac{Q_+}{Q_-}$ for $c \in \C \setminus \{0\}$.
\end{itemize}
\label{Saito}
\end{theorem}

\noindent
\textbf{Proof}. 
First assume that the condition (2) holds. We may assume that 
$\det M(\theta_1,\ldots,\theta_\ell)=\displaystyle \frac{Q_+}{Q_-}$. Note that 
this implies that $\theta_1,\ldots,\theta_\ell$ are independent over $S$. 
By definition, $Q_- \theta_i$ is a regular vector field for each $i$ and 
$\det M(Q_-\theta_1,\ldots,Q_-\theta_\ell)=Q_-^{\ell-1}Q_+$. 
Take $\theta \in D\Omega(\A,m)$. Since $Q_+Q_-^\ell \theta \in Q_+Q_-^{\ell-1} \Der_{\C}(S)$, 
there exist polynomials $f_1,\ldots,f_\ell \in S$ such that 
$Q_+Q_-^\ell \theta=\sum_{i=1}^\ell f_i Q_- \theta_i$. Hence 
it is sufficient to prove that $Q_+Q_-^{\ell-1}|f_i$ for $i=1,2,\ldots,\ell$. 
Let us compute the following (where $ \partial x=\partial_{x_1}\wedge \cdots \wedge \partial_{x_\ell})$:
\begin{eqnarray*}
Q_+Q_-^\ell \theta_1 \wedge \cdots \wedge \theta_{i-1} \wedge 
\theta \wedge \theta_{i+1} \wedge \cdots \wedge \theta_\ell 
&=& g_i Q_+^2Q_-^{\ell-1} \partial x\ (g_i \in S) \\
&=& \theta_1 \wedge \cdots \wedge (\sum_{i=1}^\ell f_i Q_- \theta_i)
\wedge \cdots \wedge \theta_\ell\\
&=& f_iQ_+  \partial x.
\end{eqnarray*}
Hence $f_i=Q_+Q_-^{\ell-1}g_i$ and $\{\theta_1,\ldots,\theta_\ell\}$ forms a basis for $D\Omega(\A,m)$. 

Next assume that the condition (1) holds. Put 
$\det M(\theta_1,\ldots,\theta_\ell)=f \displaystyle \frac{Q_+}{Q_-}$ for 
$f \in S \setminus \{0\}$ (by Lemma \ref{wedge}). Let us prove that 
$f \in \C \setminus \{0\}$. Take 
$H_0 \in \A$ and put $|m(H_0)|=m_0$. We may choose an orthonormal basis 
$\{x_1,\ldots,x_\ell\}$ for $V^*$ such that 
$\alpha_{H_0}=x_1$. First assume that 
$Q_+=x_1^{m_0} Q_+'$ with $x_1 \nmid Q_+'$. 
Then 
$Q_+ \partial_{x_1}, 
Q_+' \partial_{x_2},\ldots,Q_+' \partial_{x_\ell} 
\in D\Omega(\A,m)$. Since $\{\theta_1,\ldots,\theta_\ell\}$ forms a basis, we have 
$$
Q_+(Q_+')^{\ell-1}\partial x \in S \cdot f \displaystyle \frac{Q_+}{Q_-}\partial x,
$$
or equivalently,
$$
Q_- (Q_+')^{\ell-1} =f g\ \mbox{for some}\ g \in S.
$$
Since $H_0$ runs over all hyperplanes in $\A_+$, it follows that 
$f|Q_-$. Note that $f$ and $Q_+$ have no common factors. 
Next assume that 
$Q_-=x_1^{m_0} Q_-'$ with $x_1 \nmid Q_-'$. 
Then 
$ \displaystyle \frac{Q_+}{x_1^{m_0}} \partial_{x_1}, 
Q_+ \partial_{x_2},\ldots,Q_+ \partial_{x_\ell}
\in D\Omega(\A,m)$. Since $\{\theta_1,\ldots,\theta_\ell\}$ forms a basis, we have 
$$
\displaystyle \frac{Q_+^{\ell}}{x_1^{m_0}}\partial x \in S \cdot f \displaystyle \frac{Q_+}{Q_-}\partial x,
$$
or equivalently,
$$
Q_-' Q_+^{\ell-1} =f g'\ \mbox{for some}\ g' \in S.
$$
Since $H_0$ runs over all hyperplanes in $\A_-$, it follows that $f$ is a non-zero scalar.
\owari
\medskip

\begin{cor}
The homogeneous elements $\theta_1,\ldots,\theta_\ell \in D\Omega(\A,m)$ form a basis for 
$D\Omega(\A,m)$ if and only if the following two conditions are satisfied:
\begin{itemize}
\item[(1)] $\theta_1,\ldots,\theta_\ell$ are independent over $S$.
\item[(2)] $\sum_{i=1}^\ell \deg (\theta_i)=|m|:=\deg Q_+-\deg Q_-$.
\end{itemize}
\label{Saito2}
\end{cor}

\noindent
\textbf{Proof}. By Theorem \ref{Saito} it is obvious that 
(1) and (2) are satisfied if $\{\theta_1,\ldots,\theta_\ell\}$ forms a basis. 
Assume that (1) and (2) hold. If we put 
$\det M(\theta_1,\ldots,\theta_\ell)=f \displaystyle \frac{Q_+}{Q_-}$ for some 
$f \in S \setminus \{0\}$ (because of the condition (1)), then $\deg (f)=0$ by the condition (2). 
Hence Theorem \ref{Saito} completes the proof. \owari
\medskip

\begin{example}
Let $\A$ be the Coxeter arrangement of type $A_3$, i.e., the arrangement defined by 
$\prod_{1 \leq i<j \leq 4} (x_i-x_j)=0$ 
for an orthonormal basis $\{x_1,x_2,x_3,x_4\}$. Define the multiplicity $m_1$ on $\A$ by 
\[
m_1(H_{ij}):=
\left\{
\begin{array}{rl}
1 & \mbox{if}\ (i,j)=(1,2),\ (1,3)\ \mbox{or}\ (2,3),\\
-1 & \mbox{if}\ (i,j)=(1,4),\ (2,4)\ \mbox{or}\ (3,4),
\end{array}
\right.
\]
where $H_{ij}:=\{x_i-x_j=0\}$. In other words, $(\A,m_1)$ is defined by the rational 
function $\prod_{1 \le i<j\le 3} (x_i-x_j)/\prod_{i=1,2,3}(x_i-x_4)$. Then it is easy to check that 
the following four derivations are contained in $D\Omega(\A,m_1)$:
\begin{eqnarray*}
\theta_1:&=&\partial_{x_1}+\partial_{x_2}+\partial_{x_3},\\
\theta_2:&=&x_1\partial_{x_1}+x_2\partial_{x_2}+x_3\partial_{x_3},\\
\theta_3:&=&\partial_{x_4},\\
\theta_4:&=&\displaystyle \frac{\partial_{x_1}}{x_1-x_4}+\displaystyle \frac{\partial_{x_2}}{x_2-x_4}+
\displaystyle \frac{\partial_{x_3}}{x_3-x_4}-
(\displaystyle \frac{1}{x_1-x_4}+
\displaystyle \frac{1}{x_2-x_4}+
\displaystyle \frac{1}{x_3-x_4})\partial_{x_4}.
\end{eqnarray*}
Then Corollary \ref{Saito2} shows that $(\A,m_1)$ is free with a basis 
$\{\theta_1,\theta_2,\theta_3,\theta_4\}$ and $\exp(\A,m_1)=(0,1,0,-1)$. 
On the other hand, the multi-braid arrangement $(\A,m_2)$ 
of type $A_3$ defined by 
$(x_1-x_2)(x_2-x_3)(x_1-x_4)/(x_1-x_3)(x_2-x_4)(x_3-x_4)=0$ 
is not free (see Theorem \ref{multibraid}). 

The (non-)freeness of these mult-braid arrangements and the exponents have been 
expected by the result in \cite{ANN}, for 
the multiplicity $m_1$ corresponds to the \textit{bicolor-eliminable graph} and the exponents above 
to the degrees of that graph, but 
$m_2$ does not. 
We will consider the multi-braid arrangement of this type 
in Theorem \ref{multibraid} again. 
\label{A3}
\end{example}

Now we give the duality of logarithmic modules, which is one of the most important 
results in the arrangement theory. 

\begin{theorem}
For a multiarrangement $(\A,m)$, define 
the multiplicity $-m$ on $\A$ by $(-m)(H):=-m(H)$ for $H \in \A$. 
Equivalently, $(\A,-m)$ is defined by the rational function 
$\displaystyle \frac{Q_-}{Q_+}$. 
Then $D\Omega(\A,m)$ and $D\Omega(\A,-m)$ are $S$-dual modules.
\label{duality}
\end{theorem}

\noindent
\textbf{Proof}. 
For a basis $\{x_1,\ldots,x_\ell\}$ for $V^*$, define an $S$-bilinear map 
$\langle\ ,\ \rangle :D\Omega(\A,m) \times D\Omega(\A,-m) \rightarrow S$ by 
$\langle \theta,\omega \rangle := \sum_{i=1}^\ell \theta(x_i) \varphi(x_i)$. It is easy to check that
this pairing is independent of a choice of a basis. 
Let us prove that this pairing is a non-degenerate perfect 
pairing.
First let us prove that the image of $\langle\ ,\ \rangle$ is contained in $S$. 
Put $\theta=\sum_{i=1}^\ell f_i\partial_{x_i}/Q_-$ and 
$\omega=\sum_{i=1}^\ell g_i\partial_{x_i}/Q_+$ for $f_i,g_j \in S$. 
Assume that $Q_+=\alpha_{H_0}^{m_0} Q_+'$ with $\alpha_{H_0} \nmid Q_+'$. We prove that 
$\langle \theta,\omega \rangle=(\sum_{i=1}^\ell f_ig_i)/Q_+Q_-$ is regular along $H_0$. 
We may put $\alpha_{H_0}=x_1$ by an orthonormal change of coordinates. 
Since $\partial_{x_1} \wedge \omega$ 
is regular along $H_0$ by definition of $D\Omega(\A,-m)$, it holds that 
$g_2,\ldots,g_\ell \in S\cdot x_1^{m_0}$. Also, by definition of $D\Omega(\A,m)$, 
$Q_- \theta(x_1)=f_1$ is divisible by $x_1^{m_0}$. 
Hence $\langle \theta,\omega \rangle=
(\sum_{i=1}^\ell f_i g_i)/Q_+ Q_{-}$ is regular along 
$H_0$. 
By applying the same arguments to each $H \in \A$, we can see that 
$\langle \theta,\omega\rangle \in S$. Hence the pairing $\langle\ ,\ \rangle$ induces the $S$-homomorphism
$A:D\Omega(\A,m) \rightarrow D\Omega(\A,-m)^*:=\mbox{Hom}_S(D\Omega(\A,-m), S)$ defined by 
$A(\theta)(\omega):=\langle \theta,\omega\rangle$. Because of the symmetry, to complete the proof, it suffices to show that 
$A$ is an isomorphism. First let us prove that $A$ is injective. 
Assume that $A(\theta)=0$. Let us put $\theta=\sum_{i=1}^\ell f_i \partial_{x_i}/Q_-$. 
Noting that $Q_- \partial_{x_i} \in D\Omega(\A,-m)$ for all $i$, we can see that 
$A(\theta)(Q_- \partial_{x_i})=f_i=0$, which implies that $\theta=0$. Next let us prove that 
$A$ is surjective. Take $\varphi \in D\Omega(\A,-m)^*$ and define 
$\displaystyle \frac{\Der_{\C}(S)}{Q_-} \ni \overline{\varphi}
:=\displaystyle \frac{1}{Q_-}
\sum_{i=1}^\ell \varphi(Q_- \partial_{x_i}) \partial_{x_i}$. We prove that 
$\overline{\varphi} \in D\Omega(\A,m)$. First assume that 
$Q_+=x_1^{m_0} Q_+'$ with $x_1 \nmid Q_+'$. By definition, 
$\displaystyle \frac{Q_-}{x_1^{m_0}} \partial_{x_1} \in D\Omega(\A,-m)$. 
Hence $\overline{\varphi}(x_1)=\displaystyle \frac{1}{Q_-} \varphi(x_1^{m_0} \cdot 
\displaystyle \frac{Q_-}{x_1^{m_0}} \partial_{x_1}) \in 
S \cdot \displaystyle \frac{x_1^{m_0}}{Q_-}$. Next assume that 
$Q_-=x_1^{m_0} Q_-'$ with 
$x_1 \nmid Q_-'$ and prove that 
$\partial_{x_1} \wedge \overline{\varphi}$ is regular along $H_0$. 
It suffices to show that 
$\displaystyle \frac{1}{Q_-} \varphi(Q_- \partial_{x_i})$ is regular along 
$H_0$ for $2 \le i \le \ell$. By definition, 
$\displaystyle \frac{Q_-}{x_1^{m_0}} \partial_{x_i}
 = Q_-' \partial_{x_i} 
\in D\Omega(\A,-m)$ for $2 \le i \le\ell$. 
Hence 
$\displaystyle \frac{1}{Q_-} \varphi(Q_- \partial_{x_i})
=\displaystyle \frac{1}{Q_-} 
\varphi(x_1^{m_0} Q_-' \partial_{x_i}) 
\in S \cdot \displaystyle \frac{1}{Q_-'}$, which shows that 
$\overline{\varphi} \in D\Omega(\A,m)$. For $Q_- \partial_{x_i} \in D\Omega(\A,-m)$, 
$A( \overline{\varphi})(Q_- \partial_{x_i})=
\langle \overline{\varphi}, Q_- \partial_{x_i} \rangle=\varphi(Q_- \partial_{x_i})$. 
Hence $A(\overline{\varphi})=\varphi$ on the module 
$\sum_{i=1}^\ell S \cdot Q_- \partial_{x_i} \subset D\Omega(\A,-m)$. For a general 
$\omega = \sum_{i=1}^\ell g_i \partial_{x_i}/Q_+ \in D\Omega(\A,-m)$, it follows that 
$Q_+Q_-A(\overline{\varphi})(\omega)=
A(\overline{\varphi})(Q_+Q_-\omega)=A(\overline{\varphi})(\sum_{i=1}^\ell g_i Q_- \partial_{x_i}) 
=\varphi(\sum_{i=1}^\ell g_i Q_- \partial_{x_i})=
\varphi(Q_+Q_- \omega)=Q_+Q_-\varphi(\omega)$. Hence 
$A(\overline{\varphi})=\varphi$. \owari 
\medskip

Theorem \ref{duality} is a generalization of the well-known duality between 
$D(\A,m)$ and $\Omega^1(\A,m)$. Hence we can also obtain the following result.

\begin{cor}
For a multiarrangement $(\A,m)$, the module $D\Omega(\A,m)$ is a reflexive 
$S$-module.
\label{reflexive}
\end{cor}

\noindent
\textbf{Proof}. 
Immediate from Theorem \ref{duality}. \owari
\medskip


The following is a direct consequence of 
Corollary \ref{reflexive} and a general results on reflexive modules. 

\begin{cor}
If $(\A,m)$ is a multiarrangement in a two-dimensional vector space, then 
$(\A,m)$ is free.
\label{2multi}
\end{cor}

\begin{rem}
We can prove some other results on $D\Omega(\A,m)$, e.g., 
the decomposition of a logarithmic module of an arrangement which admits a direct sum 
decomposition (generalization of Lemma 1.4 in \cite{ATW1}), or 
the jumping behavior of a basis of free a multiarrangement and its deletion (generalization of 
Theorem 0.4 
in \cite{ATW2}). We do not give proofs since they are very easy and not used in this article. However, 
there are still a lot of results for previous multiarrangements which have not yet proved 
for generalized multiarrangements defined in this article. For instance, the preservation of 
the freeness under the localization or addition-deletion theorems have not been known yet. 
\end{rem}

\begin{example}
Note that the relation between the freeness of $(\A_+,m_+)$, $(\A_-,-m_-)$ and 
$(\A,m)$ is not clear. For example, let $\A$ be the arrangement defined by 
$(x_1-x_2)(x_2-x_3)/(x_2-x_4)$ for an orthonormal basis $\{x_1,\ldots,x_4\}$ for 
$V \simeq \R^4$. Then it is obvious that 
$(\A_+,m_+)$ and $(\A_-,-m_-)$ are both free, but we can show that 
$(\A,m)$ is not free (for example, see Theorem \ref{multibraid}). Thus 
the freeness of both $(\A_+,m_+)$ and $(\A_-,-m_-)$ does not imply that 
of $(\A,m)$. 

The natural question is whether the freeness of $(\A,m)$ imply those of 
$(\A_+,m_+)$ and $(\A_-,-m_-)$ or not. For the braid arrangement, this is true 
again by Theorem \ref{multibraid}.
\end{example}


\section{Application to Coxeter multiarrangements}
In this section we prove Theorem \ref{mainmain} and the duality of 
Coxeter multiarrangements. For that purpose, 
we apply the theory above to the freeness of Coxeter multiarrangements as done for 
quasi-constant multiplicities in \cite{AY2}, and also done for 
specific multiplicities on the braid arrangements in \cite{ANN}. 
Throughout this section, let $V$ be an $\ell$-dimensional Euclidean space, 
$W \subset \mbox{GL}(V)$ a finite Coxeter group with the Coxeter number $h$ and  
$\A=\A(W)$ the Coxeter arrangement corresponding to $W$. Consider an action of $W$ onto 
$V, V^*, S:=\mbox{Sym}^*(V^*) \otimes_\R \C,\ \Der_\C(S)$ and $\Omega^1_V$. 
Then it is known that there exists a $W$-invariant 
inner product $I^*:V^* \times V^* \rightarrow \R$ (see \cite{Sa4} for example). 
Let us fix this 
inner product in this section. Then we can consider the action of $W$ onto 
$D\Omega(\A,m)$ since the action of $W$ onto $D(\A,m)$ and $\Omega(\A,m)$ can be identified 
through the $W$-invariant inner product $I^*$. 

Now, let $S^W$ be 
the $W$-invariant subring of $S$. Then Chevalley's theorem in \cite{Ch} shows that  
$S^W= \C[P_1,\ldots,P_\ell]$ with basic invariants $P_1,\ldots,P_\ell$. It is known that 
$P_1,\ldots,P_\ell$ are algebraically independent over $\C$ and 
$2= \deg P_1 < \deg P_2 \le \ldots \le \deg P_{\ell-1} < \deg P_\ell=h$. Hence the rational 
derivation $D:= \displaystyle \frac{\partial}{\partial P_\ell}$ is unique up to scalars, and 
called the \textit{primitive derivation}. The primitive derivation plays a crucial role 
in the theory of constructing the Hodge filtration and flat structure on 
$\Der(S^W)$ in \cite{Sa1}. Also, recall that 
the affine connection $\nabla:\Der_\C(S) \times \Der_\C(S) \rightarrow 
\Der_\C(S)$ defined by 
$
\nabla_\theta(\varphi):=\sum_{i=1}^\ell \theta(\varphi(x_i)) \partial_{x_i}$ for $\theta, \varphi \in 
\Der_\C(S)$. 
Then 
applying $D$ through $\nabla$, 
several results on free Coxeter multiarrangements have been 
obtained in \cite{ST2}, \cite{T4}, \cite{T5}, \cite{Y0} and \cite{AY2}. 
We give the most recent version of these results 
through the logarithmic 
module $D\Omega(\A,m)$ as follows:

\begin{theorem}
Let $\A$ be a Coxeter arrangement with the Coxeter number $h$, 
$m:\A \rightarrow \{+1,0,-1\}$, $\theta_E:=\sum_{i=1}^\ell x_i \partial_{x_i}$  
the Euler derivation and 
$k \in \Z$. 
Define the map $\Phi_k:D\Omega(\A,m)(-kh) \rightarrow D\Omega(\A,2k+m)$ by 
$$
\Phi_k(\theta):=\nabla_\theta \nabla_D^{-k} \theta_E,  
$$
where $\theta \in D\Omega(\A,m)$. Then $\Phi_k$ 
is an $S$-module isomorphism. 
\label{Coxeter}
\end{theorem}

Note that Theorem \ref{Coxeter} 
immediately shows Theorem \ref{mainmain}. 
For the proof of Theorem \ref{Coxeter} we need the following, which 
is the dual version of Theorem 10 in \cite{AY2}.

\begin{theorem}
Let $\A$ be a Coxeter arrangement and $m:\A \rightarrow \{0,1\}$ a quasi-constant multiplicity. 
Define the morphism $\Phi_k:\Omega^1(\A,m)(-kh) \rightarrow D(\A,2k-m)\ (k \in \Z_{>0})$ by 
$$
\Phi_k(\omega):=\nabla_{I^*(\omega)} \nabla_D^{-k} \theta_E, 
$$
where $\omega \in \Omega^1(\A,m)$. 
Then $\Phi_k$ is an $S$-module isomorphism.
\label{differential}
\end{theorem}

\noindent
\textbf{Proof}. First note that the definition 
of $\Phi_k$ is independent of a choice of coordinates since 
it is defined without coordinates. 
Through the identification of $\Omega^1_V$ and $\Der_\C(S)$ we often write 
$\nabla_\omega$ (resp. $\omega(f))$ instead of $\nabla_{I^*(\omega)}$ (resp. 
$I^*(\omega)(f))$ for $f \in S$. 
We prove Theorem \ref{differential} in three steps.

\noindent
\textbf{Step 1}. Well-definedness of $\Phi_k$. First, we prove that the image of $\Phi_k$ 
is a regular vector field. Let $Q_-$ be 
the defining polynomial of $(\A,m)$ and take a hyperplane 
$H_0 \in \A$ such that $m(H_0)=1$. It suffices to show that 
$\nabla_\omega \nabla_D^{-k} \theta_E$ is regular along $H_0$ for $\omega \in \Omega^1(\A,m)$. 
We choose an orthonormal coordinates $\{x_1,\ldots,x_\ell\}$ in such a way that 
$\alpha_{H_0}=x_1$ and prove that $\nabla_\omega \nabla_D^{-k} \theta_E(x_i)$ is regular along $H_0$ for 
$i=1,2,\ldots,\ell$. 
Recall that $E_k:=\nabla_D^{-k} \theta_E \in D(\A,2k+1)$. Also, if we put 
$E_k(\alpha_{H})=\alpha_H^{2k+1} g_H\ (H \in \A,\ g_H \in S)$ then 
$\alpha_H \nmid g_H$ (see \cite{T4}, \cite{Y0} and \cite{AY2}). 
If we put $E_k(x_1)=g x_1^{2k+1}\ (g \in S)$, then $\Phi_k(\omega)(x_1)=
\omega(E_k(x_1))=\omega(g)x_1^{2k+1}+(2k+1)x_1^{2k}\omega(x_1)g$.
Since $k>0$, it is obvious that 
$\Phi_k(\omega)(x_1)$ is regular along $H_0$. 
Next consider $\Phi_k(\omega)(x_i)$ for $i \neq 1$.  
Let $\tau_1 \in W$ be the reflection corresponding to the 
hyperplane $H_0=\{x_1=0\}$. Recall that $E_k$ is $W$-invariant (see \cite{T4} and 
\cite{Y0}) and $\{x_1,\ldots,x_\ell\}$ is an orthonormal basis for $V^*$. 
Also, since $I^*$ is $W$-invariant, the actions of $W$ onto 
$dx_i$ and $\partial_{x_i}$ are the same. Thus 
$\tau_1(E_k(x_i))=\tau_1(E_k(\tau_1(x_i)))=(\tau_1 E_k)(x_i)=E_k(x_i)$. 
Hence $E_k(x_i)$ is $\tau_1$-invariant. In other words, 
$$
E_k(x_i)=\sum_{n \ge 0} f_{in}(x_2,\ldots,x_\ell)x_1^{2n}.
$$
Thus it suffices to show that $\omega(f_{i0})$ is regular along $H_0$. However, 
this is obvious 
since $\omega(x_1) I^*(d{x_1})(f_{i0})=\omega(x_1) \partial_{x_1}(f_{i0})=0$ and $\omega \wedge d {x_1}=
\sum_{i=2}^\ell \omega(x_i) d {x_i} \wedge d{x_1}$ is regular along 
$H_0$. 
Since it is easy to see that 
$\Phi_k(\omega)(\alpha_{H_0}) \in S \cdot \alpha_{H_0}^{2k-m(H_0)}$, 
the image of $\Phi_k$ is in $D(\A,2k-m)$. 

\noindent
\textbf{Step 2}. Injectivity of $\Phi_k$. Choose $S$-independent 
elements $\omega_1,\ldots,\omega_\ell \in \Omega^1(\A,m) \simeq D\Omega(\A,-m)$. 
Obviously $Q_- \omega_1,\ldots,Q_- \omega_\ell \in \Der_{\C}(S) \simeq \Omega^1_V$. Hence Lemma 7 in 
\cite{AY2} shows that 
$$
\nabla_{Q_- \omega_1} \nabla_{D}^{-k} \theta_E,
\ldots,
\nabla_{Q_- \omega_\ell} \nabla_{D}^{-k} \theta_E
\in D(\A,2k)
$$ 
are $S$-independent. Thus 
so are $\nabla_{\omega_1} \nabla_{D}^{-k} \theta_E, \ldots,\nabla_{\omega_\ell} \nabla_{D}^{-k} \theta_E$. 
Then the same argument as in the proof of Theorem 10 in 
\cite{AY2} shows the injectivity of $\Phi_k$. 

\noindent
\textbf{Step 3}. Surjectivity of $\Phi_k$. First we prove that $\Phi_k$ is isomorphic when 
$m \equiv 1$. Since $\Omega^1(\A,1)$ is free, we can take a basis $\omega_1,\ldots,\omega_\ell$ 
for $\Omega^1(\A,1)=\Omega^1(\A)$. Then $\Phi_k(\omega_1),\ldots,\Phi_k(\omega_\ell) \in D(\A,2k-1)$ 
are $S$-independent by Steps 1 and 2. Also, $\sum_{i=1}^\ell \deg (\Phi_k(\omega_i))=
kh\ell - |\A|$, which is equal to the sum of multiplicities 
$\sum_{H \in \A} (2k-1)=|\A|(2k-1)$ of $(\A,2k-1)$ 
(see \cite{Y0}). Hence Corollary \ref{Saito2} 
shows that $\Phi_k:\Omega^1(\A) \rightarrow D(\A,2k-1)$ is an isomorphism. 

Next consider an arbitrary quasi-constant multiplicity $m$ on $\A$ and 
take $\theta \in D(\A,2k-m)$. Since 
$D(\A,2k-m) \subset D(\A,2k-1)$, the previous paragraph implies that 
there exists a differential $1$-form $\omega \in \Omega^1(\A,1)$ such that 
$\nabla_\omega \nabla_D^{-k} \theta_E =\theta$. Note that $\Omega^1(\A,1) \supset \Omega^1(\A,m)$ and 
prove that $\omega \in \Omega^1(\A,m)$. By definitions of the logarithmic modules, it suffices to 
show that 
$\omega$ is regular along $H_0$ such that $m(H_0)=0$. 
Let us choose a new orthonormal coordinate system $\{x_1,\ldots,x_\ell\}$ in such a way that 
$H_0=\{x_1=0\}$. Since $\omega \in \Omega^1(\A,1)$ it holds that 
$\omega \wedge d x_1$ is regular along $x_1=0$. In other words, if we put $\omega=\sum_{i=1}^\ell 
f_i dx_i$, then $f_2,\ldots, f_\ell$ are regular along $x_1=0$. We prove that $f_1$ is also regular along $x_1$. 
Put $\nabla_D^{-k} \theta_E(x_1)=g x_1^{2k+1}$ and $\theta(x_1)=f x_1^{2k}$. Recall that 
$x_1 \nmid g$. Then 
\begin{eqnarray*}
f x_1^{2k}&=&\theta(x_1)=(\nabla_\omega \nabla_D^{-k} \theta_E)(x_1) \\
&=&x_1^{2k+1}\sum_{i=1}^\ell f_i \partial_{x_i}(g) + 
(2k+1) g x_1^{2k} f_1.
\end{eqnarray*}
Since the pole order of $f_1$ along $x_1$ is at most one, we can see that 
$f_1$ is regular along $x_1$, which completes the proof. \owari
\medskip

\noindent
\textbf{Proof of Theorem \ref{Coxeter}}. 
When $k=0$ there is nothing to prove. Since the proof is the same, 
we only prove Theorem \ref{Coxeter} when $k>0$. In other words, we prove that 
$$
\Phi_k:D\Omega(\A,m)(-kh) \rightarrow D\Omega(\A,2k+m)=D(\A,2k+m)
$$
is an $S$-isomorphism. 
Since the same proofs as in Theorem \ref{differential} on the well-definedness and injectivity 
are valid in the setup of Theorem \ref{Coxeter}, it suffices to show that 
$\Phi_k$ is surjective. 
Take $\varphi \in D(\A,2k+m)$ and recall the definition of $m_\mu\ (\mu \in \{+,-\})$ 
in Definition \ref{DO}. Then 
$D(\A,2k+m)=D(\A,2k+m_++m_-) \subset D(\A,2k+m_-)$. Hence Theorem \ref{differential} 
shows that there exists a logarithmic $1$-form $\omega \in \Omega^1(\A,-m_-)$ such that 
$\nabla_\omega \nabla_D^{-k} \theta_E=\varphi$. 
By using this identification and a canonical isomorphism 
$I^*:\Omega^1(\A,-m_-) \simeq D\Omega(\A,m_-)$, the differential form $\omega$ can be expressed 
as $\theta/Q_- \in D\Omega(\A,m_-)$ for $\theta \in \Der_\C(S)$, see the following commutative 
diagram:

$$\xymatrix@R1pc{
D\Omega(\A,m) \ar[d]_{\Phi_k}\ \ar @{^{(}-{>}}[r] &\ D\Omega(\A,m_-)
\ar[d]_{\Phi_k
} \ar[r]^{(I^*)^{-1}} & \Omega^1(\A,-m_-) \ar[d]_{\Phi_k
} \\
D\Omega(\A,2k+m)\ \   \ar @{^{(}-{>}}[r] & \ \ D\Omega(\A,2k+m_-) \ar@{=}[r] & D\Omega(\A,2k+m_-)}$$

To complete the proof, 
it suffices to show that 
$\omega=\displaystyle \frac{\theta}{Q_-} \in D\Omega(\A,m)$, or equivalently, 
$\theta(\alpha_H) \in S \cdot \alpha_H$ for any $H \in \A_+$. 
Take $H \in \A_+$ and put $\varphi(\alpha_H)=f \alpha_H^{2k+1}\ (f \in S)$. 
Because $\nabla_{\theta/Q_-}\nabla_D^{-k} \theta_E =\varphi$, 
it holds that $Q_- f \alpha_H^{2k+1}=\theta(\nabla_D^{-k}\theta_E(\alpha_H))$. Put 
$\nabla_D^{-k} \theta_E(\alpha_H)=g \alpha_H^{2k+1}$ with $g \in S$. 
Recall that 
$\alpha_H \nmid g$. 
Therefore, 
$$
Q_- f \alpha_H^{2k+1}=\theta(g\alpha_H^{2k+1})=
\theta(g)\alpha_H^{2k+1}+(2k+1)\alpha_H^{2k}\theta(\alpha_H)g.
$$
Hence $\alpha_H \mid \theta(\alpha_H)$, which completes the proof. \owari
\medskip

\begin{example}
Let $\A$ be the Coxeter arrangement of type $B_2$ defined by $xy(x^2-y^2)=0$, where 
$\{x,y\}$ is an orthonormal coordinate system for $\R^2$. 
For the $B_2$-type, $P_1=(x_1^2+x_2^2)/2$ and $P_2=(x^4+y^4)/4$ for 
$S^W=\C[P_1,P_2]$. Moreover, 
$D=(-y \partial_x +x \partial_y)/(x^3y-xy^3)$ and 
$\nabla_D^{-1}\theta_E=(-x^5/15 + x^3y^2/3) \partial_x+(x^2y^3/3-y^5/15) \partial_y$.
As we have seen in Example \ref{B2}, the multiarrangement defined by 
$y/x(x-y)=0$ is free with a basis
$$
\theta_1=\displaystyle \frac{\partial_x}{x},\ 
\theta_2=\displaystyle \frac{x \partial_x - y \partial_y}{x-y}.
$$
Hence, for the multiarrangement $(\A,m)$ defined by $xy^3(x-y)(x+y)^2=0$, Theorem \ref{Coxeter} 
shows that $D(\A,m)=\langle \nabla_{\theta_1} \nabla_D^{-1} \theta_E,
 \nabla_{\theta_2} \nabla_D^{-1} \theta_E \rangle_S.$

\end{example}




As a corollary of Theorem \ref{Coxeter}, 
the duality theorem in \cite{AY2} is also generalized as follows:

\begin{cor}
Let $\A$ be a Coxeter arrangement with the Coxeter number $h$ and 
$m:\A \rightarrow \{+1,0,-1\}$. 
Then $D(\A,2k+m)(kh)$ and $D(\A,2k-m)(kh)$ are 
$S$-dual modules.
\label{coxeterduality}
\end{cor}

By using Definition \ref{DO} and Theorem \ref{Coxeter}, the result in \cite{ANN} can be understood in terms of 
a characterization of the freeness of $D\Omega(\A,m)$ for the braid arrangement $\A$. Let 
$m$ be a multiplicity on the braid arrangement $\A$ with 
$\mbox{Im}(m) \subset \{+1,0,-1\}$. Then the set of these multiplicities has 
a one to one correspondence with an edge-bicolored graph. Then  
Theorem 0.3 in \cite{ANN} can be extended as follows:

\begin{theorem}
Let $k \in \Z$. 
In the above notation, the multi-braid arrangement $(\A,2k+m)$ is free if and only if 
the corresponding graph with $m$ is bicolor-eliminable.
\label{multibraid}
\end{theorem}

\noindent
\textbf{Proof}. 
By Theorem \ref{duality} and result in \cite{ANN} it suffices to show the case when 
$k=0$. Then Theorem \ref{Coxeter} shows that 
$(\A,m)$ is free if and only if $(\A,2k+m)$ is free for any $k \in \Z_{>0}$. Hence 
Theorem 0.3 in \cite{ANN} completes the proof. \owari
\medskip

\begin{rem}
The results in \cite{Sa1}, 
\cite{Sa3} and \cite{Sa4} (see also \cite{T4} and \cite{Y0})  
give us the following commutative diagram:

$$\xymatrix@R1pc{
\Omega^1_{S^W}  \ar[d]_{I^*} \ar[r]^{\nabla_D} & \Omega^1(\A)^W \ar[d]_{I^*} \\
\Der_\R(S)^W \ar[r]^{\nabla_D} &  \Der_{S^W} }$$

Here both rows are $T(:=\R[P_1,\ldots,P_{\ell-1}])$-isomorphisms and 
both columns are $S^W$-isomorphisms. It is known that 
$\Der_\R(S)^W \simeq D(\A)^W$. Also, it is proved by Terao in \cite{T6} and 
\cite{T7} that the second row can be extended to the left as 
a sequence of shift isomorphisms:
$$
\cdots \rightarrow D(\A, 2k+1)^W 
\stackrel{\nabla_D}{\rightarrow} D(\A,2k-1)^W 
\stackrel{\nabla_D}{\rightarrow} D(\A,2k-3)^W \rightarrow \cdots .
$$
Combining this sequence with the above diagram and identification $I^*$, it 
seems natural that the image $\nabla_D(D(\A)^W) \simeq \Der_{S^W} \simeq 
I^*(\Omega^1(\A)^W)$ might be regarded as $D(\A,-1)^W$, which also supports 
Definition \ref{DO}.
\end{rem}

 \vspace{5mm}

\end{document}